\documentclass[a4paper,10pt]{article}
\usepackage{amscd}
\usepackage{bbding}
\usepackage{amsfonts}
\usepackage{graphicx,amssymb,amsthm}
\usepackage{mathrsfs}
\usepackage{latexsym,amsmath,cite,enumerate}
\usepackage{float}
\usepackage{leftidx}
\usepackage{manfnt}
\usepackage{stmaryrd}
\usepackage{cases}
\usepackage{bbm}
\usepackage{tabls}
\usepackage{multirow}
\usepackage[dvipdfm,
           pdfstartview=FitH,
           colorlinks=true,
             pdfborder=001,
           citecolor=blue]{hyperref}

\allowdisplaybreaks[2]
\DeclareFontEncoding{OT2}{}{}
\DeclareTextSymbol{\cyrsftsn}{OT2}{126}
\DeclareTextSymbol{\textnumero}{OT2}{125}

\theoremstyle{definition}
\newtheorem{theorem}{Theorem}[section]
\newtheorem{lemma}{Lemma}[section]
\newtheorem{corollary}{Corollary}[section]

\newtheorem{definition}{Definition}[section]
\newtheorem{remark}{Remark}[section]
\newtheorem{example}{Example}[section]

\begin{document}

\title{{\LARGE \textbf{A New System of Global Fractional-order Interval Implicit Projection Neural Networks}}\thanks{This work was supported by
the National Natural Science Foundation of China (11471230, 11671282).}}
\author{Zeng-bao Wu$^{a,b}$, Jin-dong Li$^{b,c}$, Nan-jing Huang$^b$\thanks{Corresponding author. E-mail addresses: njhuang@scu.edu.cn; nanjinghuang@hotmail.com} \\
%EndAName
{\small\it $^a$Department of Mathematics, Luoyang Normal University, Luoyang, Henan  471022, P.R. China}\\
{\small\it $^b$Department of Mathematics, Sichuan University, Chengdu, Sichuan 610064, P.R. China}\\
{\small\it $^c$College of Management Science, Chengdu University of Technology, Chengdu, Sichuan 610059, P.R. China}}
\date{ }
\maketitle
\begin{flushleft}
\hrulefill\newline
\end{flushleft}
\textbf{Abstract.} The purpose of this paper is to introduce and investigate a new system of global fractional-order interval implicit projection neural networks. An existence and uniqueness theorem of the equilibrium point for the system of global fractional-order interval implicit projection neural networks is obtained under some suitable assumptions. Moreover, Mittag-Leffler stability for the system of global fractional-order interval implicit projection neural networks is also proved. Finally, two numerical examples are given to illustrate the validity of our results. \newline
\
\newline
\textbf{Key Words and Phrases:} Interval implicit projection neural networks, fractional-order calculus, equilibrium point, Mittag-Leffler stability.

\begin{flushleft}
\hrulefill
\end{flushleft}

\section{Introduction}
\noindent \setcounter{equation}{0}
This paper deals with a new system of global fractional-order interval implicit projection neural networks (FIIPNN) in $R^n \times R^m$ as the following form:
\begin{equation}
\left\{
\begin{array}{l}
\leftidx{_0^C}D{_t^\alpha}x(t) =P_{K_1(x(t))}[x(t) -\rho \left(Ax(t)+A^\ast y(t)\right)-\rho a]-x(t),\ t\geq 0, \\
x(0)=x_0=(x_{10}, x_{20}, \ldots,x_{n0})^\top, \\
\leftidx{_0^C}D{_t^\alpha}y(t) =P_{K_2(y(t))}[y(t)-\lambda \left(By(t)+B^\ast x(t)\right)-\lambda b]-y(t),\ t\geq 0, \\
y(0) =y_0=(y_{10}, y_{20}, \ldots,y_{m0})^\top,
\end{array}
\right.   \label{FIIPNN}
\end{equation}
where $\alpha \in (0,1)$, $\leftidx{_0^C}D{_t^\alpha}$ is the Caputo fractional derivative, $K_1: R^n \rightarrow 2^{R^n}$ and $K_2: R^m \rightarrow 2^{R^m}$ are two point to set mappings with nonempty, closed and convex values, $P_{K_1(x(t))}$ and $P_{K_2(y(t))}$ are two implicit projection operators, $\rho >0$ and $\lambda >0$ are two constants, $a=(a_1, a_2, \ldots,a_n)^\top \in R^n$ and $b=(b_1, b_2, \ldots,b_m)^\top\in R^m$ are two vectors, and
\begin{equation*}
\left\{
\begin{array}{l}
A\in A_{I}=\left\{\left. \left(a_{ij}\right)_{n\times n}\right|\underline{A}\leq A\leq \overline{A},\; \text{i.e.},\; \underline{a}_{ij}\leq a_{ij}\leq \overline{a}_{ij} \right\} , \\
A^{\ast }\in A_{I}^{\ast }=\left\{\left.\left(a_{ij}^{\ast }\right)_{n\times m}\right| \underline{A^{\ast }}\leq A^{\ast }\leq \overline{A^{\ast }},\; \text{i.e.},\;
\underline{a^{\ast }}_{ij}\leq a_{ij}^{\ast }\leq \overline{a^{\ast }}_{ij} \right\} , \\
B\in B_{I}=\left\{ \left.\left( b_{ij}\right) _{m\times m}\ \right|\underline{B}\leq B\leq \overline{B},\; \text{i.e.},\; \underline{b}_{ij}\leq b_{ij}\leq
\overline{b}_{ij}  \right\} , \\
B^{\ast }\in B_{I}^{\ast }=\left\{\left. \left( b_{ij}^{\ast }\right) _{m\times n}\right| \underline{B^{\ast }}\leq B^{\ast }\leq \overline{B^{\ast }},\; \text{i.e.},\;
\underline{b^{\ast }}_{ij}\leq b_{ij}^{\ast }\leq \overline{b^{\ast }}_{ij} \right\}.
\end{array}
\right.
\end{equation*}

Some special cases of (\ref{FIIPNN}) are as follows.
\begin{itemize}
\item[(i)]  If $K_1(u)\equiv K_1$ and $K_2(v)\equiv K_2$ for all $(u,v)\in R^n\times R^m$, here $K_1\subset R^n$ and $K_2\subset R^m$ are two nonempty, closed and convex subsets, then (\ref{FIIPNN}) reduces to the following problem:
\begin{equation}
\left\{
\begin{array}{l}
\leftidx{_0^C}D{_t^\alpha}x(t) =P_{K_1}[x(t) -\rho \left(Ax(t)+A^\ast y(t)\right)-\rho a]-x(t),\ t\geq 0, \\
x(0)=x_0=(x_{10}, x_{20}, \ldots,x_{n0})^\top, \\
\leftidx{_0^C}D{_t^\alpha}y(t) =P_{K_2}[y(t)-\lambda \left(By(t)+B^\ast x(t)\right)-\lambda b]-y(t),\ t\geq 0, \\
y(0) =y_0=(y_{10}, y_{20}, \ldots,y_{m0})^\top,
\end{array}
\right.   \label{e1.2}
\end{equation}
which is the system of fractional-order generalized projection neural networks introduced and studied by Wu et al. \cite{WZH2}.

\item[(ii)] If $n=m$, $\underline{A}=A=\overline{A}=\underline{B}=B=\overline{B}$, $\underline{A^\ast}=A^\ast =\overline{A^\ast}=0$, $\underline{B^\ast}=B^\ast =\overline{B^\ast}=0$, $a=b$, $\rho=\lambda$, $x_0=y_0$ and $K_1(u)=K_2(u)\equiv K$ for all $u\in R^n$, here $K\subset R^n$ is a nonempty, closed and convex subset, then (\ref{FIIPNN}) reduces to the following problem:
\begin{equation}
\left \{
\begin{array}{l}
\leftidx{_0^C}D{_t^\alpha}x(t)=P_{K}[ x(t) -\rho Ax(t)-\rho a ] -x(t), \quad t\geq0,  \\
x_i(0) =x_{i0}, \quad i=1,2,\cdots,n,
\end{array}
\right.   \label{e1.3}
\end{equation}
which is the global projection dynamical systems with fractional-order introduced and investigated by Wu and Zou \cite{WZ}.

\item[(iii)] If $\alpha =1$, $n=m$, $\underline{A}=A=\overline{A}=\underline{B}=B=\overline{B}$, $\underline{A^\ast}=A^\ast =\overline{A^\ast}=0$, $\underline{B^\ast}=B^\ast =\overline{B^\ast}=0$, $a=b$, $\rho=\lambda$, $x_0=y_0$ and $K_1(u)=K_2(u)$ for all $u\in R^n$, then (\ref{FIIPNN}) reduces to the implicit projected dynamical systems considered by Noor et al. \cite{NNK}.
\end{itemize}

We remark that for suitable choices of $K_1$, $K_2$, $A$, $A^*$, $B$, $B^*$ and $\alpha$ in the formulation of
\eqref{FIIPNN}, one can obtain many problems of the fractional-order projection neural networks (dynamical systems) and the implicit projection neural networks (dynamical systems) investigated in recent literature.

We note that the projection neural networks (dynamical systems) have been used to solve constrained optimization problems, variational inequality problems, complementarity problems, dynamic traffic network and  interregional commodity movements and so on \cite{ASA,JB,DN,D,FBMT,FSB,DG,WWZ,WZLX,WWZ2,WZH1,XV,Xia,Z,ZLHS,ZS,ZWZS}. On the other hand, fractional order systems have also been a hot research topic due to their application for control theory, mechanics and physics, viscoelasticity materials, biology, electrical circuits, neural networks and so on (see,  for example, \cite{AHJ,LZ,OK,IP,SMM,TB}). Nevertheless, in some practical world, it is necessary to consider some complex systems such as model involving \eqref{FIIPNN}. In the following, we will present an example, which comes from \cite{FBMT}.

\begin{example}\label{e1.1}
A network tatonnement model was introduced by Friesz et al. \cite{FBMT} to investigate dynamics of network adjustments. In particular, we study a simple network model involving 5 arcs and 6 nodes (see, Figure \ref{arc_node}), which has a origin (node 1) and a destination (node 4) with three paths. Path $p_1$ is composed of arcs $a_1$ and $a_4$, path $p_2$ is composed of arcs $a_2$, $a_3$ and $a_4$, path $p_3$ is composed of arcs $a_2$ and $a_5$. Here, we follow the notations used in \cite{WZH2}.
\begin{figure}[H]
\begin{center}
\includegraphics[width=2.5in]{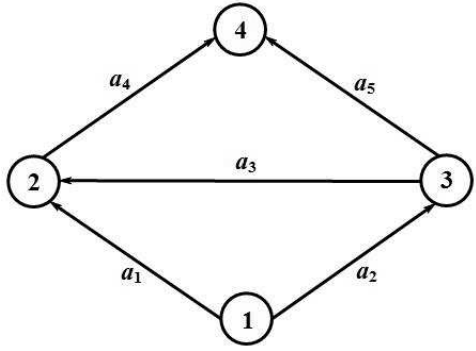}
\end{center}
\caption{5-arcs, 4-nodes traffic network.}
\label{arc_node}
\end{figure}
Applying the network tatonnement model presented by Friesz et al. \cite{FBMT}, we have
\begin{equation}
\left \{
\begin{array}{l}
\frac{dh_{p_1}(t)}{dt} = \kappa_1 \{P_{K_{1,1}}[h_{p_1}(t) - \rho (C_{a_1}(h_p(t))+C_{a_4}(h_p(t))-u_{14}(t))] -h_{p_1}(t)\}, \\
\frac{dh_{p_2}(t)}{dt} = \kappa_2 \{P_{K_{1,2}}[h_{p_2}(t) - \rho (C_{a_2}(h_p(t)) + C_{a_3}(h_p(t)) + C_{a_4}(h_p(t)) - u_{14}(t))] -h_{p_2}(t)\}, \\
\frac{dh_{p_3}(t)}{dt} = \kappa_3 \{P_{K_{1,3}}[h_{p_3}(t) - \rho (C_{a_2}(h_p(t)) +C_{a_5}(h_p(t)) - u_{14}(t))] - h_{p_3}(t)\}, \\
\frac{du_{14}(t)}{dt} = \eta_1 \{P_{K_{2,1}}[u_{14}(t) -\lambda (h_{p_1}(t)+h_{p_2}(t)+h_{p_3}(t)-T_{14}(u_{14}(t)))] -u_{14}(t)\},
\\
h_{p_1}(0) =h_{p_1}^0,\; h_{p_2}(0) =h_{p_2}^0,\; h_{p_3}(0) =h_{p_3}^0,\; u_{14}(0) =u_{14}^0,\;  t\geq 0,
\end{array}
\right.  \label{STM}
\end{equation}
where $\rho$, $\lambda$, $\kappa_i \,(i=1,2,3)$, $\eta_1$ are positive constants,  $h_p=(h_{p_1}, h_{p_2}, h_{p_3})^\top$, $K_{1,i} \,(i=1,2,3)$, $K_{2,1}$ denote the feasible constraints with fixed lower and upper bounds for  flows $h_{p_i}\,(i=1,2,3)$ and cost $u_{14}$, respectively, that is,
$$K_{1,i}=\left\{h_{p_i} \in R | \, c_{1,i} \leq h_{p_i} \leq c_{2,i} \right\} \, (i=1,2,3), $$
$$K_{2,1}=\left\{u_{14} \in R | \, d_1 \leq u_{14} \leq d_2 \right\}. $$
With the adjustment of flows and cost, it is difficult to maintain the same fixed bounds for constraint sets. Thus, it is reasonable to assume that $c_{1,i}$ and $c_{2,i}$ are dependent on flows, $d_1$ and $d_2$ are dependent on travel cost, that is, $c_{1,i}$ and $c_{2,i}$ are functions of $h_p$, $d_1$ and $d_2$ are functions of $u_{14}$. Therefore, the constraint sets can be rewritten as follows:
$$K_{1,i}(h_p)=\left\{\varrho_i \in R | \, c_{1,i}(h_p) \leq \varrho_i \leq c_{2,i}(h_p)\right\} \, (i=1,2,3), $$ $$K_{2,1}(u_{14})=\left\{\nu \in R | \, d_1(u_{14}) \leq \nu \leq d_2(u_{14})\right\}. $$

We assume that the cost functions of flow on arc $a_m$ can be written as follows
\begin{eqnarray*}
C_{a_m}(h_p(t))=l_m \cdot \sum_{i=1}^3\chi_{a_{m}p_{i}}h_{p_i}(t), \quad m=1,2,3,4,5
\end{eqnarray*}
and the travel demand function can  be written as
\begin{eqnarray*}
T_{14}(u_{14}(t))=r \cdot u_{14}(t),
\end{eqnarray*}
where $r$, $l_m \, (m=1,2,3,4,5)$ are real numbers, $\chi_{a_{m}p_{i}}=1$ if $a_{m} \in p_{i}$ and $\chi_{a_{m}p_{i}}=0$ otherwise. Unfortunately, it is difficult to determine the precise values of these coefficients in practice, whereas it is easier to give certain confidence intervals for these coefficients, namely, $\underline{r} \leq r \leq \overline{r}$ and $\underline{l}_m \leq l_m \leq \overline{l}_m$.

Moreover, as Example 1.1 of \cite{WZH2} indicates, this dynamic network has memory. We observe that, for the problem with memory, it is more appropriate to use the fractional order model rather than integer one (see, for instance, \cite{LZ,OK,TB,SMM}). By the above discussion, we know that model (\ref{STM}) can be reformed as the following fractional order form:
\begin{equation}
\left \{
\begin{array}{l}
\leftidx{_0^C}D{_t^\alpha} h_{p_1}(t)= \kappa_1 \{P_{K_{1,1}(h_p(t))}[h_{p_1}(t) - \rho ((l_1+l_4)h_{p_1}(t)+l_4 h_{p_2}(t)-u_{14}(t))] -h_{p_1}(t)\}, \\
\leftidx{_0^C}D{_t^\alpha} h_{p_2}(t) = \kappa_2 \{P_{K_{1,2}(h_p(t))}[h_{p_2}(t) - \rho (l_4 h_{p_1}(t)+\sum_{m=2}^4 l_i h_{p_2}(t)+l_2 h_{p_3}(t) - u_{14}(t))] -h_{p_2}(t)\}, \\
\leftidx{_0^C}D{_t^\alpha} h_{p_3}(t) = \kappa_3 \{P_{K_{1,3}(h_p(t))}[h_{p_3}(t) - \rho (l_2 h_{p_2}(t)+ (l_2+l_5)h_{p_3}(t) -u_{14}(t))] -h_{p_3}(t)\}, \\
\leftidx{_0^C}D{_t^\alpha} u_{14}(t) = \eta_1 \{P_{K_{2,1}(u_{14}(t))}[u_{14}(t) -\lambda (h_{p_1}(t) + h_{p_2}(t) + h_{p_3}(t)-r u_{14}(t))] -u_{14}(t)\},
\\
h_{p_1}(0) =h_{p_1}^0,\; h_{p_2}(0) =h_{p_2}^0,\; h_{p_3}(0) =h_{p_3}^0,\; u_{14}(0) =u_{14}^0,\;  t\geq 0,
\end{array}
\right.  \label{GSTM}
\end{equation}
where $\underline{l}_m \leq l_m \leq \overline{l}_m$ and $\underline{r} \leq r \leq \overline{r}$.

Clearly, if $\kappa_i=1$ $(i=1,2,3)$ and $\eta_1=1$, then model (\ref{GSTM}) is a form of (\ref{FIIPNN}).
\end{example}

It is worth to mention that FIIPNN (\ref{FIIPNN}) is fascinating and important both as its equilibrium behavior is depicted by the quasi variational inequality (QVI for short), and also because the equilibrium point set of FIIPNN (\ref{FIIPNN}) coincides with the solution set to a QVI problem. It is well known that QVI problem is an important generalization of the variational inequality problems (see, for example,  \cite{BC,GLT,Mosco}).  Furthermore, we note that FIIPNN (\ref{FIIPNN}) obtains the desired properties of both the fractional-order system and the QVI within the same framework. Consequently, it is meaningful to investigate the equilibrium point of FIIPNN (\ref{FIIPNN}) and the stability for FIIPNN (\ref{FIIPNN}).  The main purpose of this paper is to give some new conditions to guarantee the existence and uniqueness of the equilibrium point for FIIPNN (\ref{FIIPNN}),  and the new stability result for FIIPNN (\ref{FIIPNN}) which improves some known stability results in \cite{WZ,WZH2}.

The outline of this paper is as follows. Some definitions and known results are presented in Section 2.  The existence and uniqueness concerned with the equilibrium point for FIIPNN (\ref{FIIPNN}) and the stability results in connection with the FIIPNN (\ref{FIIPNN}) are showed in Section 3.  Finally, two numerical examples to demonstrate the main conclusions are given in Section 4.

\section{Preliminaries}
\noindent \setcounter{equation}{0}

In this section, we first recall some known definitions and facts.

Following the definitions of \cite{IP,Die,AHJ}, the Riemann-Liouville fractional integral with order $\alpha>0$ is described as
\begin{equation*}
I^{\alpha}_{t_0}x(t)=\frac{1}{\Gamma(\alpha)}\int_{t_0}^t(t-\tau)^{\alpha-1}x(\tau)d\tau, \quad t>t_0,
\end{equation*}
where $\Gamma(\cdot)$ is the gamma function, and the Caputo fractional derivative with order $\alpha \in (0,1)$ is described as
\begin{equation*}
\leftidx{_{ t_0}^C}D{_t^\alpha}x(t)=\frac{1}{\Gamma(1-\alpha)}\int_{t_0}^t(t-\tau)^{-\alpha}x'(\tau)d\tau,\quad t>t_0.
\end{equation*}
Moreover, the Mittag-Leffler function with two parameters $\alpha>0$ and $\beta>0$ is defined by
$$E_{\alpha,\beta}(z)=\sum_{k=0}^\infty\frac{z^k}{\Gamma(\alpha k+\beta)},\quad \alpha>0, \; \beta>0, \; z \in C.$$
For $\beta=1$, the one-parameter Mittag-Leffler function is shown as
$$E_\alpha(z):=E_{\alpha,1}(z)=\sum_{k=0}^\infty\frac{z^k}{\Gamma(\alpha k+1)}, \quad \alpha>0, \; z \in C.$$
In particular, $E_1(z)=e^z.$

\begin{definition}
\label{def2.4} Assume that $K:R^n\to 2^{R^n}$ is a point to set mapping with nonempty, closed and convex values. For any given $x\in  R^n$, the implicit projection operator $P_{K(x)}:R^n \rightarrow K(x)$ is described as
\begin{equation*}
P_{K(x)}[y] =\arg\!\!\!\min_{z\in K(x)}\Vert y-z\Vert, \quad y \in R^n.
\end{equation*}
\end{definition}
\begin{remark}
\label{uxK} In many applications \cite{BC,GLT,Mosco,CP82}, the point to set mapping $K(x)$ can be given by the following form:
$$K(x)=u(x)+K,$$
where $K \subset R^n$ is a closed convex set, $u(x)$ is a point to point mapping and the addition of a point $v$ and a set $K$ is defined by $v+K=\{v+w| \, w \in K\}$. In this case, the following relation holds
$$P_{K(x)}[y]=P_{u(x)+K}[y]=u(x)+P_K[y-u(x)], \quad \forall \; x, \, y \in R^n.$$
\end{remark}

\begin{definition}
\label{def2.5} A vector $(x^\ast,y^\ast) \in R^n \times R^m$ is called an equilibrium point of  (\ref{FIIPNN}) if, for each $A \in A_I$, $A^\ast \in A^\ast_I$, $B \in B_I$ and $B^\ast \in B^\ast_I$, the vector $(x^\ast,y^\ast)$ satisfies the following relations:
\begin{equation}
\left \{
\begin{array}{l}
P_{K_1\left(x^\ast\right)}\left[x^\ast -\rho \left(Ax^\ast+A^\ast y^\ast\right)-\rho a\right]=x^\ast, \\
P_{K_2\left(y^\ast\right)}\left[y^\ast-\lambda \left(By^\ast+B^\ast x^\ast\right)-\lambda b\right]=y^\ast.
\end{array}
\right.  \notag
\end{equation}

\end{definition}

\begin{lemma}
\label{nonexpan}\cite[Corollary 2.4]{DG} If $K$ is a convex closed subset of a Hilbert space $H$, then the projection $P_K$ is non-expansive, i.e.,
\begin{equation*}
\left\|P_K[u]-P_K[v]\right\| \leq \|u-v\|, \quad \forall u,v \in H.
\end{equation*}
\end{lemma}

\begin{lemma}
\label{e-unique}\cite[Remark 3.8]{LCP} Let $G=[t_0,+\infty)\times \Omega$. Assume that $g:G\rightarrow R^n$ is continuous such that it fulfils the locally  Lipschitz condition with respect to the second variable, where $\Omega \subset R^n$ is a domain. Then there exists a unique solution $x(t)$ of the following initial-value problem
\begin{equation}\label{fDg}
\left\{
\begin{array}{l}
\leftidx{_{ t_0}^C}D{_t^\alpha}x(t)=g(t,x), \; \alpha \in (0,1], \\
x(t_0)=x_0.
\end{array}
\right.
\end{equation}
\end{lemma}

\begin{definition}
(Mittag-Leffler Stability \cite{LCP}) If $x^\ast$ is an equilibrium point of (\ref{fDg}), then the solution of (\ref{fDg}) is called Mittag-Leffler stable if there exist two constants $\lambda>0$ and $b>0$ such that
\begin{equation*}
\left\|x(t)- x^\ast\right\| \leq [m\left(x(t_0)-x^\ast\right) E_\alpha\left(-\lambda (t-t_0)^\alpha\right)]^b,
\end{equation*}
where $m(0)=0$, $m(x)\geq 0$, and $m(x)$ is locally Lipschitz on $x \in R^n$.
\end{definition}

\begin{lemma}
\label{Dsgn}\cite[Theorem 2]{ZYW} If $x(t) \in C^1([0,+\infty),R)$, then
\begin{equation*}
\leftidx{_0^C}D{_t^\alpha}|x\left(t^+\right)| \leq sgn(x(t))\leftidx{_0^C}D{_t^\alpha}x(t) \; \text{(holding almost everywhere)},
\end{equation*}
where $0<\alpha<1$ and $x(t^+):=\lim_{s\rightarrow t^+}x(s)$.
\end{lemma}

\begin{lemma}
\label{MLS}\cite[Theorem 1]{ZYW} For $t_0=0$, let $V(t,x(t)):[0,+\infty) \times \Omega \rightarrow R$ be a continuous function satisfying the locally Lipschitzian condition with respect to the second variable such that
$$\alpha_1\|x(t)\|^a \leq V(t,x(t)) \leq \alpha_2\|x(t)\|^{ab},$$
$$\leftidx{_0^C}D{_t^\beta}V\left(t^+,x\left(t^+\right)\right) \leq -\alpha_3\|x(t)\|^{ab} \; \text{(holding almost everywhere)},$$
where $\dot{V}(t,x(t))$ is piecewise continuous, $\lim_{s\rightarrow t^+}\dot{V}(s,x(s))$ exists for any $t \in [0,+\infty)$, $\Omega \subset R^n$ is a domain containing the origin and $V\left(t^+,x\left(t^+\right)\right):=\lim_{s\rightarrow t^+}V(s,x(s))$, $t\geq0$, $\beta \in (0,1)$, $\alpha_i$ $(i=1,2,3)$, $a$ and $b$ are positive constants. Then system (\ref{fDg}) is Mittag-Leffler stable at the equilibrium point $x^\ast=0$. Moreover, if all the assumptions are satisfied globally on $R^n$, then system (\ref{fDg}) is globally  Mittag-Leffler stable at the equilibrium point $x^\ast=0$.
\end{lemma}

\section{Main results}
\noindent \setcounter{equation}{0}

From now on we make the following assumptions:
\begin{enumerate}[(\mbox{A}$_1$)]
  \item For any $x=(x_1,x_2,\ldots,x_n)^\top \in R^n$ and $y=(y_1,y_2,\ldots,y_m)^\top \in R^m$, $K_1(x)$ and $K_2(y)$ are assumed to be as follows
  $$K_1(x)=u^1(x)+K_1,\quad K_2(y)=u^2(y)+K_2,$$
 where
  $$u^1(x)=\left(u^1_1(x),u^1_2(x),\ldots,u^1_n(x)\right)^\top, \quad u^1_i(x)=\sum_{j=1}^n h_{ij}x_j, \quad i=1,2,\ldots,n,$$
  $$u^2(y)=\left(u^2_1(y),u^2_2(y),\ldots,u^2_m(y)\right)^\top, \quad u^2_j(y)=\sum_{i=1}^m l_{ji}y_i, \quad j=1,2,\ldots,m, $$
  and
  $$
   K_1=\left\{x \in R^n| \, c_{1,i}\leq x_i\leq c_{2,i}, \; i=1,2,\ldots,n\right\},\quad
 K_2=\left\{y \in R^m| \, d_{1,j}\leq y_j\leq d_{2,j}, \; j=1,2,\ldots,m\right\},
$$
here $h_{ij}$, $l_{ji}$, $c_{1,i}$, $c_{2,i}$, $d_{1,j}$ and $d_{2,j}$ are all constants;
  \item $1 - \rho \overline{a}_{ii} \geq h_{ii},\ i=1,2,\cdots,n$;
  \item $1 - \lambda \overline{b}_{jj} \geq l_{jj}, \ j=1,2, \cdots, m$;
  \item There exist constants $\mu_i>0 \; (i=1,2,\ldots,n)$ and $\tau_j>0 \; (j=1,2,\ldots,m)$ such that
  \begin{eqnarray*}\label{A2}
 \left \{
\begin{array}{l}
0 < 1- \rho \underline{a}_{ii}-h_{ii}+|h_{ii}|+ \sum \limits_{j=1,j\neq i}^n\frac{\mu_j}{\mu_i} \left(\widetilde{a}_{ji}+|h_{ji}|\right)+ \sum \limits_{j=1}^{m}\frac{\tau_j}{\mu_i} \lambda \widetilde{b^{\ast}}_{ji} < 1 , \ i=1,2,\cdots,n, \\
0 < 1- \lambda \underline{b}_{jj}-l_{jj}+|l_{jj}|+ \sum \limits_{i=1,i\neq j}^m\frac{\tau_i}{\tau_j} \left(\widetilde{b}_{ij}+|l_{ij}|\right)+ \sum \limits_{i=1}^n\frac{\mu_i}{\tau_j} \rho \widetilde{a^{\ast}}_{ij} < 1 , \ j=1,2, \cdots, m,
\end{array}
\right.
  \end{eqnarray*}
where
$$
\widetilde{a}_{ji}=\max\left\{\left| \rho\underline{a}_{ji}+h_{ji}\right|,  \left|\rho\overline{a}_{ji}+h_{ji}\right|\right\}, \quad \widetilde{a^\ast}_{ij}=\max\left\{\left|\underline{a^\ast}_{ij}\right|, \left|\overline{a^\ast}_{ij}\right|\right\}
$$
and
$$
\widetilde{b}_{ij}=\max\left\{\left|\lambda\underline{b}_{ij}+l_{ij}\right|,\left|\lambda\overline{b}_{ij}+l_{ij}\right|\right\}, \quad
\widetilde{b^\ast}_{ji}=\max\left\{\left|\underline{b^\ast}_{ji}\right|, \left|\overline{b^\ast}_{ji}\right|\right\}.
$$
\end{enumerate}

Clearly, under assumption (A$_1$),  an equivalent formulation of (\ref{FIIPNN}) can be rewritten as follows
\begin{equation}
\left\{
\begin{array}{l}
\leftidx{_0^C}D{_t^\alpha}x_i(t) =P_{K_{1,i}(x(t))}\left[x_i(t) -\rho \left(\sum\limits_{j=1}^n a_{ij}x_j(t)+\sum\limits_{j=1}^m a_{ij}^\ast y_j(t)\right) -\rho a_i \right] -x_i(t),\ t\geq 0, \\
x_i\left( 0\right) =x_{i0},\ i=1,2,\cdots ,n, \\
\leftidx{_0^C}D{_t^\alpha}y_j(t) =P_{K_{2,j}(y(t))}\left[y_j(t) -\lambda \left(\sum\limits_{i=1}^m b_{ji}y_i(t)+\sum\limits_{i=1}^n b_{ji}^\ast x_i(t)\right) -\lambda b_j \right] -y_j(t),\ t\geq 0, \\
y_j\left( 0\right) =y_{j0},\ j=1,2,\cdots ,m, \\
\end{array}
\right.   \label{FIIPNN1}
\end{equation}
where
$\underline{a}_{ij}\leq a_{ij}\leq \overline{a}_{ij}, \;
\underline{a^\ast }_{ij}\leq a^\ast_{ij}\leq \overline{a^\ast}_{ij}, \;
\underline{b}_{ji}\leq b_{ji}\leq \overline{b}_{ji}, \;
\underline{b^\ast}_{ji}\leq b_{ji}^\ast \leq \overline{b^\ast}_{ji}$,
\begin{equation}\label{Kix}
K_{1,i}(x(t))=u^1_i(x(t))+K_{1,i}=\sum_{j=1}^n h_{ij}x_j(t)+K_{1,i}, \quad i=1,2,\ldots,n,
\end{equation}
and
\begin{equation}\label{Kiy}
K_{2,j}(y(t))=u^2_j(y(t))+K_{2,j}=\sum_{i=1}^m l_{ji}y_i(t)+K_{2,j},  \quad j=1,2,\ldots,m,
\end{equation}
with
\begin{equation}\label{K1x}
K_{1,i}=\{x_i \in R| \, c_{1,i}\leq x_i\leq c_{2,i}\}, \quad i=1,2,\ldots,n,
\end{equation}
and
\begin{equation}\label{K2y}
K_{2,j}=\{y_j \in R| \, d_{1,j}\leq y_j\leq d_{2,j}\}, \quad j=1,2,\ldots,m.
\end{equation}

\subsection{Existence and uniqueness of the equilibrium point}

This subsection will present an existence and uniqueness theorem concerned with the equilibrium point for (\ref{FIIPNN}).

\begin{theorem}
\label{Exist} Assume that all assumptions (A$_1$)-(A$_4$) are satisfied. Then FIIPNN (\ref{FIIPNN}) has a unique equilibrium point for each $A \in A_I$, $A^\ast \in A^\ast_I$, $B \in B_I$ and $B^\ast \in B^\ast_I$.

\noindent\textbf{Proof.}\hspace{0.2cm} For any given $A \in A_I$, $A^\ast \in A^\ast_I$, $B \in B_I$ and $B^\ast \in B^\ast_I$, let $T_{\rho i}: R^n \times R^m \rightarrow R$ be given by
\begin{eqnarray}
 T_{\rho i}(x,y)=\mu_i P_{K^{\mu}_{1,i}(x)}\left[ \frac{x_i}{\mu_i} -\rho \left(\sum\limits_{j=1}^{n}a_{ij}\frac{x_j}{\mu_j}+\sum\limits_{j=1}^{m}a_{ij}^\ast \frac{y_j}{\tau_j}\right)-\rho a_i\right], \quad \forall \ (x,y) \in R^n \times R^m, \label{Toi}
\end{eqnarray}
and $T_\rho(x,y)$ be given by
\begin{eqnarray}
T_\rho(x,y)=\left(T_{\rho 1}(x,y),T_{\rho 2}(x,y),\cdots,T_{\rho n}(x,y)\right)^\top, \quad \forall \ (x,y) \in R^n \times R^m,\label{To}
\end{eqnarray}
where
$$K^{\mu}_{1,i}(x)=u^{1,\mu}_i(x)+K_{1,i},$$
$$u^{1,\mu}_i(x)=\sum_{j=1}^n h_{ij}\frac{x_j}{\mu_j}, \quad \mu=\left(\mu_1,\mu_2,\ldots,\mu_n\right)^\top$$
and $K_{1,i}$ is defined by (\ref{K1x}).

Let
$$
\Vert x\Vert =\sum_{i=1}^n\vert x_{i}\vert ,\quad \forall x=(x_{1},x_{2},\ldots,x_{n})^\top \in R^n.
$$
For any given vectors $\left(x^1,y^1\right)$ and $\left(x^2,y^2\right)$ in $R^{n} \times R^{m}$, by (\ref{Toi}), (\ref{To}) and Remark \ref{uxK}, one has
\begin{eqnarray}
&& \left\| T_\rho\left(x^2,y^2\right)-T_\rho\left(x^1,y^1\right) \right\|  \notag \\
&=& \sum_{i=1}^n \left| T_{\rho i}\left(x^2,y^2\right)-T_{\rho i}\left(x^1,y^1\right) \right|  \notag \\
&=& \sum_{i=1}^n \mu_i \left|P_{K^\mu_{1,i}\left(x^2\right)}\left[ \frac{x_i^2}{\mu_i} -\rho \left(\sum\limits_{j=1}^n a_{ij}\frac{x_j^2}{\mu_j} + \sum\limits_{j=1}^m a_{ij}^\ast \frac{y_j^2}{\tau_j}\right)-\rho a_i \right]  \notag  \right.\\
&&\left. - P_{K^\mu_{1,i}\left(x^1\right)}\left[ \frac{x_i^1}{\mu_i} -\rho \left(\sum\limits_{j=1}^n a_{ij}\frac{x_j^1}{\mu_j} + \sum\limits_{j=1}^m a_{ij}^\ast \frac{y_j^1}{\tau_j}\right)-\rho a_i\right] \right|  \notag \\
&=& \sum_{i=1}^n\mu_i \left| u^{1,\mu}_i\left(x^2\right) + P_{K_{1,i}}\left[ \frac{x_i^2}{\mu_i} -\rho \left(\sum\limits_{j=1}^n a_{ij}\frac{x_j^2}{\mu_j} + \sum\limits_{j=1}^m a_{ij}^\ast \frac{y_j^2}{\tau_j}\right)-\rho a_i-u^{1,\mu}_i\left(x^2\right)\right] \right.  \notag \\
&&\left. - u^{1,\mu}_i\left(x^1\right) - P_{K_{1,i}}\left[\frac{x_i^1}{\mu_i} -\rho \left(\sum\limits_{j=1}^n a_{ij}\frac{x_j^1}{\mu_j} + \sum\limits_{j=1}^m a_{ij}^\ast \frac{y_j^1}{\tau_j}\right)-\rho a_i- u^{1,\mu}_i\left(x^1\right)\right] \right|  \notag \\
&=& \sum_{i=1}^n\mu_i \left| \sum_{j=1}^n h_{ij}\frac{x^2_j-x^1_j}{\mu_j} + P_{K_{1,i}}\left[ \frac{x_i^2}{\mu_i} - \sum\limits_{j=1}^n \left(\rho a_{ij}+h_{ij}\right)\frac{x_j^2}{\mu_j} - \sum\limits_{j=1}^{m}\rho a_{ij}^\ast \frac{y_j^2}{\tau_j}-\rho a_i\right] \right. \notag \\
&&\left. -P_{K_{1,i}}\left[ \frac{x_i^1}{\mu_i} -\sum\limits_{j=1}^n \left(\rho a_{ij}+h_{ij}\right)\frac{x_j^1}{\mu_j} - \sum\limits_{j=1}^m \rho a_{ij}^\ast \frac{y_j^1}{\tau_j}-\rho a_i \right] \right|  \notag \\
&\leq& \sum_{i=1}^n\mu_i \left\{ \sum_{j=1}^n |h_{ij}|\frac{|x^2_j-x^1_j|}{\mu_j} + \left|P_{K_{1,i}}\left[ \frac{x_i^2}{\mu_i} - \sum\limits_{j=1}^n \left(\rho a_{ij}+h_{ij}\right)\frac{x_j^2}{\mu_j} - \sum\limits_{j=1}^m \rho a_{ij}^\ast \frac{y_j^2}{\tau_j}-\rho a_i\right] \right.\right. \notag \\
&&\left.\left. -P_{K_{1,i}}\left[ \frac{x_i^1}{\mu_i} -\sum\limits_{j=1}^n \left(\rho a_{ij}+h_{ij}\right)\frac{x_j^1}{\mu_j} - \sum\limits_{j=1}^{m}\rho a_{ij}^\ast \frac{y_j^1}{\tau_j}-\rho a_i \right]\right| \right\}. \label{TO1}
\end{eqnarray}
By assumptions (A$_2$) and (A$_4$),  it follows from Lemma \ref{nonexpan} that
\begin{eqnarray}
&&\left|P_{K_{1,i}}\left[ \frac{x_i^2}{\mu_i} - \sum\limits_{j=1}^n \left(\rho a_{ij}+h_{ij}\right)\frac{x_j^2}{\mu_j} - \sum\limits_{j=1}^m \rho a_{ij}^\ast \frac{y_j^2}{\tau_j}-\rho a_i\right] \right.\notag \\
&&\left. -P_{K_{1,i}}\left[ \frac{x_i^1}{\mu_i} -\sum\limits_{j=1}^n \left(\rho a_{ij}+h_{ij}\right)\frac{x_j^1}{\mu_j} - \sum\limits_{j=1}^{m}\rho a_{ij}^\ast \frac{y_j^1}{\tau_j}-\rho a_i \right]\right|  \notag \\
&\leq& \left|\left[ \frac{x_i^2}{\mu_i} - \sum\limits_{j=1}^n \left(\rho a_{ij}+h_{ij}\right)\frac{x_j^2}{\mu_j} - \sum\limits_{j=1}^m \rho a_{ij}^\ast \frac{y_j^2}{\tau_j}-\rho a_i\right] \right. \notag \\
&&\left. -\left[ \frac{x_i^1}{\mu_i} -\sum\limits_{j=1}^n \left(\rho a_{ij}+h_{ij}\right)\frac{x_j^1}{\mu_j} - \sum\limits_{j=1}^m \rho a_{ij}^{\ast}\frac{y_j^1}{\tau_j}-\rho a_i \right]\right|   \notag \\
&=& \left| \frac{x_i^2-x_i^1}{\mu_i} - \sum\limits_{j=1}^n \left(\rho a_{ij}+h_{ij}\right)\frac{x_j^2-x_j^1}{\mu_j} - \sum\limits_{j=1}^m \rho a_{ij}^\ast \frac{y_j^2-y_j^1}{\tau_j} \right|  \notag \\
&\leq&  \sum\limits_{j=1,j\neq i}^n\left|\rho a_{ij}+h_{ij}\right|\frac{\left|x_j^2-x_j^1\right|}{\mu_j} + \left|1-\rho a_{ii}-h_{ii}\right|\frac{\left|x_i^2-x_i^1\right|}{\mu_i}  + \sum\limits_{j=1}^m\rho \left|a_{ij}^\ast \right|\frac{\left|y_j^2-y_j^1\right|}{\tau_j} \notag \\
&\leq& \sum_{j=1,j\neq i}^n \widetilde{a}_{ij} \frac{\left|x^2_j-x^1_j\right|}{\mu_j}  + \left(1-\rho \underline{a}_{ii} - h_{ii}\right)\frac{\left|x_i^2-x_i^1\right|}{\mu_i} + \sum\limits_{j=1}^m \rho \widetilde{a^\ast}_{ij}\frac{\left|y_j^2-y_j^1\right|}{\tau_j}. \label{TO2}
\end{eqnarray}
In light of (\ref{TO1}) and (\ref{TO2}), we have
\begin{eqnarray}
&& \left\| T_\rho\left(x^2,y^2\right)-T_\rho\left(x^1,y^1\right) \right\|  \notag \\
&\leq& \sum_{i=1}^n \left\{ \sum_{j=1,j\neq i}^n \frac{\mu_i}{\mu_j}\left(|h_{ij}| + \widetilde{a}_{ij}\right) \left|x^2_j-x^1_j\right|  + \left(|h_{ii}|+1-\rho \underline{a}_{ii} - h_{ii}\right)\left|x_i^2-x_i^1\right| \right\} \notag \\
&& + \sum_{i=1}^n\sum\limits_{j=1}^m \frac{\mu_i}{\tau_j} \rho \widetilde{a^\ast}_{ij}\left|y_j^2-y_j^1\right| \notag \\
&=& \sum_{i=1}^n \left\{ \sum_{j=1,j\neq i}^n \frac{\mu_j}{\mu_i}\left(|h_{ji}| + \widetilde{a}_{ji}\right) + |h_{ii}|+1-\rho \underline{a}_{ii} - h_{ii} \right\} \left|x_i^2-x_i^1\right| \notag \\
&& + \sum_{j=1}^m\sum\limits_{i=1}^n \frac{\mu_i}{\tau_j} \rho \widetilde{a^\ast}_{ij}\left|y_j^2-y_j^1\right|. \label{TO}
\end{eqnarray}

Moreover, let $T_{\lambda j}:R^n \times R^m \rightarrow R$ be given by
\begin{eqnarray}
 T_{\lambda j}(x,y)=\tau_j P_{K^\tau_{2,j}(y)}\left[ \frac{y_j}{\tau_j} -\lambda \left(\sum\limits_{i=1}^{m}b_{ji}\frac{y_i}{\tau_i}+\sum\limits_{i=1}^{n}b_{ji}^\ast\frac{x_i}{\mu_i}\right)-\lambda b_j\right], \ \forall \ (x,y) \in R^n \times R^m, \label{Tai}
\end{eqnarray}
and $T_\lambda(x,y)$ be given by
\begin{eqnarray}
T_\lambda(x,y)=\left(T_{\lambda 1}(x,y),T_{\lambda 2}(x,y),\cdots,T_{\lambda m}(x,y)\right)^\top, \ \forall \ (x,y) \in R^n \times R^m,\label{Ta}
\end{eqnarray}
where
$$K^\tau_{2,j}(y)=u^{2,\tau}_j(y)+K_{2,j},$$
$$u^{2,\tau}_j(y)=\sum_{i=1}^m l_{ji}\frac{y_i}{\tau_i}, \quad \tau=\left(\tau_1,\tau_2,\ldots,\tau_m\right)^\top$$
and $K_{2,j}$ is defined by (\ref{K2y}).
Then as same as the proof of (\ref{TO}), by assumptions (A$_3$) and (A$_4$), it follows from (\ref{Tai}), (\ref{Ta}), Remark \ref{uxK} and Lemma \ref{nonexpan} that
\begin{eqnarray}
&& \left\| T_\lambda\left(x^2,y^2\right)-T_\lambda\left(x^1,y^1\right) \right\| = \sum\limits_{j=1}^{m} \left| T_{\lambda j}\left(x^2,y^2\right)-T_{\lambda j}\left(x^1,y^1\right) \right| \notag \\
&\leq& \sum_{j=1}^m \left\{ \sum_{i=1,i\neq j}^m \frac{\tau_i}{\tau_j}\left(|l_{ij}| + \widetilde{b}_{ij}\right) + |l_{jj}|+1-\lambda \underline{b}_{jj} - l_{jj} \right\} \left|y_j^2-y_j^1\right| \notag \\
&& + \sum_{i=1}^n\sum\limits_{j=1}^m \frac{\tau_j}{\mu_i} \lambda \widetilde{b^\ast}_{ji}\left|x_i^2-x_i^1\right|. \label{TA}
\end{eqnarray}
Combining (\ref{TO}) and (\ref{TA}), one has
\begin{eqnarray}
&&\left\| T_\rho\left(x^2,y^2\right)-T_\rho\left(x^1,y^1\right) \right\| + \left\| T_\lambda\left(x^2,y^2\right)-T_\lambda\left(x^1,y^1\right) \right\|  \notag \\
&\leq& \sum_{i=1}^n \left\{ \sum_{j=1,j\neq i}^n \frac{\mu_j}{\mu_i}\left(|h_{ji}| + \widetilde{a}_{ji}\right) + |h_{ii}|+1-\rho \underline{a}_{ii} - h_{ii} \right\} \left|x_i^2-x_i^1\right| \notag \\
&& + \sum_{j=1}^m\sum\limits_{i=1}^n \frac{\mu_i}{\tau_j} \rho \widetilde{a^\ast}_{ij}\left|y_j^2-y_j^1\right|
+ \sum_{i=1}^n\sum\limits_{j=1}^m \frac{\tau_j}{\mu_i} \lambda \widetilde{b^\ast}_{ji}\left|x_i^2-x_i^1\right| \notag \\
&&+ \sum_{j=1}^m \left\{ \sum_{i=1,i\neq j}^m \frac{\tau_i}{\tau_j}\left(|l_{ij}| + \widetilde{b}_{ij}\right) + |l_{jj}|+1-\lambda \underline{b}_{jj} - l_{jj} \right\} \left|y_j^2-y_j^1\right| \notag \\
&=& \sum \limits_{i=1}^n \xi_i \left|x_i^2-x_i^1\right| +\sum\limits_{j=1}^m \zeta_j \left|y_j^2-y_j^1\right|  \notag \\
&\leq& \kappa\left(\left\|x^2-x^1\right\|+\left\|y^2-y^1\right\|\right), \label{TOA}
\end{eqnarray}
where
\begin{equation}\label{xi}
\xi_i=\sum_{j=1,j\neq i}^n \frac{\mu_j}{\mu_i}\left(|h_{ji}| + \widetilde{a}_{ji}\right)+ \sum\limits_{j=1}^m \frac{\tau_j}{\mu_i} \lambda \widetilde{b^\ast}_{ji} + |h_{ii}|+1-\rho \underline{a}_{ii} - h_{ii}, \; i=1,2,\ldots,n,
\end{equation}
and
\begin{equation}\label{zeta}
\zeta_j=\sum_{i=1,i\neq j}^m \frac{\tau_i}{\tau_j}\left(|l_{ij}| + \widetilde{b}_{ij}\right)+\sum\limits_{i=1}^n \frac{\mu_i}{\tau_j} \rho \widetilde{a^\ast}_{ij} + |l_{jj}|+1-\lambda \underline{b}_{jj} - l_{jj}, \; j=1,2,\ldots,m,
\end{equation}
with
$$
\kappa=\max\left\{\max\limits_{1\leq i\leq n} \xi_{i},\max\limits_{1\leq j\leq m} \zeta_{j}\right\}.
$$
It follows from assumption (A$_4$) that $0<\kappa<1$.

Let $T_{\rho\lambda}:R^n \times R^m\rightarrow R^n \times R^m$ define by $T_{\rho\lambda}(x,y)=\left(T_\rho(x,y),T_\lambda(x,y)\right)$. For any $x \in R^n$ and $y\in R^m$, let $\|(x,y)\|_1=\|x\|+\|y\|$. Then, it is well known that $\left(R^n \times R^m,\|\cdot\|_1\right)$ is a Banach space. Now (\ref{TOA}) implies that
\begin{eqnarray}
&&\left\Vert T_{\rho\lambda}\left(x^2,y^2\right)- T_{\rho\lambda}\left(x^1,y^1\right)\right\Vert  \notag \\
&=&\left\Vert \left(T_\rho\left(x^2,y^2\right),T_\lambda\left(x^2,y^2\right)\right) - \left(T_\rho\left(x^1,y^1\right),T_\lambda\left(x^1,y^1\right)\right)\right\Vert   \notag \\
&=&\left\Vert T_\rho\left(x^2,y^2\right)-T_\rho\left(x^1,y^1\right) \right\| + \left\|T_\lambda\left(x^2,y^2\right)-T_\lambda\left(x^1,y^1\right) \right\Vert  \notag \\
&\leq&\kappa\left(\left\|x^2-x^1\right\|+\left\|y^2-y^1\right\|\right) \notag \\
&=&\kappa\left\Vert\left(x^2,y^2\right)-\left(x^1,y^1\right)\right\Vert.  \label{T}
\end{eqnarray}
Thus, \eqref{T} shows that $T_{\rho\lambda}$ is contractive and therefore there exists a unique $\left(u^\ast,v^\ast\right) \in R^n \times R^m$ such that $T_{\rho\lambda}\left(u^\ast,v^\ast\right)=\left(u^\ast,v^\ast\right)$, that is,
\begin{equation}
\left \{
\begin{array}{l}
\mu_{i}P_{K^\mu_{1,i}(u^\ast)}\left[ \frac{u_i^{\ast}}{\mu_i}-\rho \left(\sum\limits_{j=1}^n a_{ij} \frac{u_j^\ast}{\mu_{j}} + \sum
\limits_{j=1}^m a^\ast_{ij} \frac{v_j^\ast}{\tau_j}\right)-\rho a_i \right]=u_i^\ast, \ i=1, 2,\ldots , n, \\
\tau_{j}P_{K^\tau_{2,j}(v^\ast)}\left[ \frac{v_j^\ast}{\tau_j}-\lambda \left(\sum\limits_{i=1}^m b_{ji} \frac{v_i^\ast}{\tau_i} + \sum
\limits_{i=1}^n b^\ast_{ji}\frac{u_i^\ast}{\mu_i}\right)-\lambda b_j \right]=v_j^\ast, \ j=1, 2,\ldots , m,
\end{array}
\right.  \notag
\end{equation}
for each $A \in A_I$, $A^\ast \in A^\ast_I$, $B \in B_I$ and $B^\ast \in B^\ast_I$. Let $x_i^\ast=\frac{u_i^\ast}{\mu_i} \ (i=1,2,\cdots,n)$ and $y_j^\ast=\frac{v_j^\ast}{\tau_j} \ (j=1,2,\cdots,m)$. Then, it is easy to see that
\begin{equation}
\left \{
\begin{array}{l}
P_{K_1\left(x^\ast\right)}\left[ x^{\ast}-\rho \left(A x^\ast + A^\ast y^\ast \right)-\rho a \right] = x^{\ast}, \\
P_{K_2\left(y^\ast\right)}\left[ y^{\ast}-\lambda \left( B y^\ast+ B^\ast x^\ast \right)-\lambda b \right]
= y^{\ast}
\end{array}
\right.  \notag
\end{equation}
and so the proof is complete. \hfill $\square$
\end{theorem}

\begin{remark}
We note that Theorem \ref{Exist} is a generalization of Theorem 3.1 in \cite{WZH2}.
\end{remark}

\begin{remark}
\label{SolExist} From the definition of $T_{\rho\lambda}$ and (\ref{T}), it is easy to check that
$$
\left\|(T_{\rho\lambda}-I)\left(x^2,y^2\right)-(T_{\rho\lambda}-I)\left(x^1,y^1\right)\right\| \leq (1+\kappa)\left\|\left(x^2,y^2\right)-\left(x^1,y^1\right)\right\|, \quad \forall \
\left(x^1,y^1\right), \left(x^2,y^2\right) \in R^n \times R^m,
$$
where $I$ is an identity mapping. Thus, we know that $T_{\rho\lambda}-I$ is a Lipschitzian mapping and so Lemma \ref{e-unique} shows that there exists a unique solution for FIIPNN (\ref{FIIPNN}) for any given $A \in A_I$, $A^\ast \in A^\ast_I$, $B \in B_I$ and $B^\ast \in B^\ast_I$.
\end{remark}

\subsection{Global Mittag-Leffler stability}
In this subsection, we will show that FIIPNN (\ref{FIIPNN}) is globally Mittag-Leffler stable under some mild conditions.

\begin{theorem}
\label{Stab} Assume that all the assumptions (A$_1$)-(A$_4$) are satisfied. Then FIIPNN (\ref{FIIPNN}) is globally Mittag-Leffler stable for each $A \in A_I$, $A^\ast \in A^\ast_I$, $B \in B_I$ and $B^\ast \in B^\ast_I$.

\noindent \textbf{Proof.}\hspace{0.2cm}  For any given $A \in A_I$, $A^\ast \in A^\ast_I$, $B \in B_I$ and $B^\ast \in B^\ast_I$. According to Remark \ref{SolExist}, we deduce that FIIPNN (\ref{FIIPNN}) has a unique solution.  Assume that
$$
\left(x^1(t),y^1(t)\right)=\left(\left(x_1^1(t),x_2^1(t),\cdots,x_n^1(t)\right)^\top, \left(y_1^1(t),y_2^1(t),\cdots,y_m^1(t)\right)^\top\right)
$$
and
$$
\left(x^{2}(t),y^{2}(t)\right)=\left(\left(x_{1}^{2}(t),x_{2}^{2}(t),\cdots,x_{n}^{2}(t)\right)^\top,
\left(y_{1}^{2}(t),y_{2}^{2}(t),\cdots,y_{m}^{2}(t)\right)^\top\right)
$$
are two solutions of FIIPNN (\ref{FIIPNN}) with different initial values
$$
\left(x^1(0),y^1(0)\right)=\left(\left(x_1^1(0),x_2^1(0),\cdots,x_n^1(0)\right)^\top,
\left(y_1^1(0),y_2^1(0),\cdots,y_m^1(0)\right)^\top\right)$$
and
$$
\left(x^2(0),y^2(0)\right)=\left(\left(x_1^2(0),x_2^2(0),\cdots,x_n^2(0)\right)^\top,
\left(y_1^2(0),y_2^2(0),\cdots,y_m^2(0)\right)^\top\right),
$$
respectively.  Let
$$e_i^1(t)=x_i^2(t)-x_i^1(t) \; (i=1,2,\cdots,n), \quad  e^1(t)=\left(e_1^1(t), e_2^1(t), \cdots, e_n^1(t)\right)^\top$$
and
$$e_j^2(t)=y_j^2(t)-y_j^1(t) \; (j=1,2,\cdots,m), \quad  e^2(t)=\left(e_1^2(t), e_2^2(t), \cdots, e_m^2(t)\right)^\top.$$
In light of assumption (A$_1$), it follows from (\ref{FIIPNN1}) and Remark \ref{uxK} that
\begin{eqnarray}
&&\leftidx{_0^C}D{_t^\alpha}e^1_i(t)  \notag\\
&=& P_{K_{1,i}\left(x^2(t)\right)}\left[x_i^2(t) -\rho \left(\sum\limits_{j=1}^n a_{ij}x_j^2(t)+\sum\limits_{j=1}^m a_{ij}^\ast y_j^2(t)\right) -\rho a_i \right]  \notag\\
&&- P_{K_{1,i}\left(x^1(t)\right)}\left[x_i^1(t) -\rho \left(\sum\limits_{j=1}^n a_{ij}x_j^1(t)+\sum\limits_{j=1}^m a_{ij}^\ast y_j^1(t)\right) -\rho a_i \right] -e^1_i(t) \notag\\
&=& u^1_i\left(x^2(t)\right)-e^1_i(t) +P_{K_{1,i}}\left[x_i^2(t) -\rho \left(\sum\limits_{j=1}^n a_{ij}x_j^2(t)+\sum\limits_{j=1}^m a_{ij}^\ast y_j^2(t)\right) -\rho a_i- u^1_i\left(x^2(t)\right)\right] \notag\\
&&   -u^1_i\left(x^1(t)\right)- P_{K_{1,i}}\left[x_i^1(t) -\rho \left(\sum\limits_{j=1}^n a_{ij}x_j^1(t)+\sum\limits_{j=1}^m a_{ij}^\ast y_j^1(t)\right) -\rho a_i-u^1_i\left(x^1(t)\right) \right] \notag\\
&=& \sum_{j=1}^n h_{ij}e^1_j(t)-e^1_i(t)+P_{K_{1,i}}\left[x_i^2(t) -\sum\limits_{j=1}^n \left(\rho a_{ij} + h_{ij}\right) x_j^2(t) - \sum\limits_{j=1}^m \rho a_{ij}^\ast y_j^2(t) -\rho a_i\right] \notag\\
&&- P_{K_{1,i}}\left[x_i^1(t) -\sum\limits_{j=1}^n \left(\rho a_{ij} + h_{ij}\right)x_j^1(t) - \sum\limits_{j=1}^m \rho a_{ij}^\ast y_j^1(t) -\rho a_i\right].  \label{Dei1}
\end{eqnarray}
Similarly, we can show that
\begin{eqnarray}
&&\leftidx{_0^C}D{_t^\alpha}e^2_j(t)
=P_{K_{2,j}}\left[y_j^2(t) -\sum\limits_{i=1}^m \left(\lambda b_{ji} + l_{ji}\right) y_i^2(t) - \sum\limits_{i=1}^n \lambda b_{ji}^\ast x_i^2(t) -\lambda b_j\right] \notag\\
&&- P_{K_{2,j}}\left[y_j^1(t) -\sum\limits_{i=1}^m \left(\lambda b_{ji} + l_{ji}\right)y_i^1(t) - \sum\limits_{i=1}^n \lambda b_{ji}^\ast x_i^1(t) -\lambda b_j\right]+ \sum_{i=1}^m l_{ji}e^2_i(t)-e^2_j(t).  \label{Dei2}
\end{eqnarray}

Let $\mu_{i}>0$ ($i=1,2,\ldots,n$), $\tau_{j}>0$  ($j=1,2,\ldots,m$),
\begin{equation}
V\left(t,\left(e^1(t),e^2(t)\right)\right)= \sum \limits_{i=1}^n \mu_i\left|e_i^1(t)\right|+\sum \limits_{j=1}^m \tau_j\left|e_j^2(t)\right| \label{Vt}
\end{equation}
and
\begin{equation}
V_1\left(t,e^1(t)\right)= \sum \limits_{i=1}^n \mu_i \left|e^1_i(t)\right|. \label{Vt1}
\end{equation}
Then, applying Lemmas \ref{nonexpan} and \ref{Dsgn} with assumptions (A$_2$) and (A$_4$), it follows from (\ref{Dei1}) and (\ref{Vt1}) that
\begin{eqnarray}
&&\leftidx{_0^C}D{_t^\alpha}V_1\left(t^+,e^1\left(t^+\right)\right) =\sum \limits_{i=1}^n\mu_i \, \leftidx{_0^C}D{_t^\alpha}\!\left|e^1_i(t^+)\right| \leq \sum\limits_{i=1}^{n}\mu_i \, sgn\left(e^1_i (t)\right) \, \leftidx{_0^C}D{_t^\alpha}e^1_i (t)  \notag \\
&=& \sum \limits_{i=1}^n \mu_i \, sgn\left(e^1_i(t)\right)\left\{P_{K_{1,i}}\left[x_i^2(t) -\sum\limits_{j=1}^n \left(\rho a_{ij} + h_{ij}\right) x_j^2(t)-\sum\limits_{j=1}^m \rho a_{ij}^\ast y_j^2(t) -\rho a_i\right] \right.\notag\\
&&\left.- P_{K_{1,i}}\left[x_i^1(t) -\sum\limits_{j=1}^n \left(\rho a_{ij} + h_{ij}\right)x_j^1(t)-\sum\limits_{j=1}^m \rho a_{ij}^\ast y_j^1(t) -\rho a_i\right]+ \sum_{j=1}^n h_{ij}e^1_j(t)-e^1_i(t)\right\} \notag \\
&\leq& \sum \limits_{i=1}^n\mu_i\left\{ \left| P_{K_{1,i}}\left[x_i^2(t) -\sum\limits_{j=1}^n \left(\rho a_{ij} + h_{ij}\right) x_j^2(t)-\sum\limits_{j=1}^m \rho a_{ij}^\ast y_j^2(t) -\rho a_i\right]   \right. \right.  \notag \\
&&\left. \left. - P_{K_{1,i}}\left[x_i^1(t) -\sum\limits_{j=1}^n \left(\rho a_{ij} + h_{ij}\right)x_j^1(t)-\sum\limits_{j=1}^m \rho a_{ij}^\ast y_j^1(t) -\rho a_i\right] \right|+ \sum_{j=1}^n \left|h_{ij}\right|\left|e^1_j(t)\right|-\left|e^1_i(t)\right| \right\} \notag \\
&\leq& \sum \limits_{i=1}^n \mu_i\left\{\left| \left[x_i^2(t) -\sum\limits_{j=1}^n \left(\rho a_{ij} + h_{ij}\right) x_j^2(t)-\sum\limits_{j=1}^m \rho a_{ij}^\ast y_j^2(t) -\rho a_i\right] \right.\right.  \notag \\
&&\left.\left. -\left[x_i^1(t) -\sum\limits_{j=1}^n \left(\rho a_{ij} + h_{ij}\right)x_j^1(t)-\sum\limits_{j=1}^m \rho a_{ij}^\ast y_j^1(t) -\rho a_i\right] \right|  + \sum_{j=1}^n \left|h_{ij}\right| \left|e^1_j(t)\right| - \left|e^1_i(t)\right|  \right\}  \notag \\
&=& \sum \limits_{i=1}^n\mu_i\left\{\left|e_i^1(t) -\sum\limits_{j=1}^n \left(\rho a_{ij} + h_{ij}\right) e_j^1(t)-\sum\limits_{j=1}^m \rho a_{ij}^\ast e_j^2(t)\right|   + \sum_{j=1}^n \left|h_{ij}\right| \left|e^1_j(t)\right| - \left|e^1_i(t)\right| \right\}  \notag \\
&=& \sum \limits_{i=1}^n\mu_i\left\{\sum_{j=1,j\neq i}^n \left|h_{ij}\right|\left|e^1_j(t)\right| +\left|e_i^1(t) -\sum\limits_{j=1}^n \left(\rho a_{ij} + h_{ij}\right) e_j^1(t)-\sum\limits_{j=1}^m \rho a_{ij}^\ast e_j^2(t)  \right| \right\}  \notag \\
&& + \sum \limits_{i=1}^n \mu_i\left(\left|h_{ii}\right|-1\right)\left|e^1_i(t)\right|  \label{DV1}
\end{eqnarray}
and
\begin{eqnarray}
&& \sum \limits_{i=1}^n\mu_i\left\{\sum_{j=1,j\neq i}^n \left|h_{ij}\right|\left|e^1_j(t)\right| +\left|e_i^1(t) -\sum\limits_{j=1}^n \left(\rho a_{ij} + h_{ij}\right) e_j^1(t) - \sum\limits_{j=1}^m \rho a_{ij}^\ast e_j^2(t) \right| \right\}  \notag \\
&\leq& \sum \limits_{i=1}^n\mu_i \left\{\sum_{j=1,j\neq i}^n \left(\left|h_{ij}\right|+\left|\rho a_{ij} + h_{ij}\right|\right) \left|e^1_j(t)\right| +\left|1-\left(\rho a_{ii}+h_{ii}\right) \right|\left|e^1_i(t) \right| +\sum\limits_{j=1}^m \rho \left|a_{ij}^\ast\right| \left|e_j^2(t)\right| \right\}  \notag \\
&\leq& \sum \limits_{i=1}^n \left\{\sum_{j=1,j\neq i}^n \mu_i\left(\left|h_{ij}\right|+\widetilde{a}_{ij}\right) \left|e^1_j(t)\right| + \mu_i\left(1-\rho \underline{a}_{ii}-h_{ii} \right)\left|e^1_i(t)\right| +\sum\limits_{j=1}^m \mu_i\rho \widetilde{a^{\ast}}_{ij} \left|e_j^2(t)\right| \right\}  \notag \\
&=& \sum \limits_{i=1}^n \left\{\sum_{j=1,j\neq i}^n \frac{\mu_j}{\mu_i} \left(\left|h_{ji}\right|+\widetilde{a}_{ji}\right) \mu_i \left|e^1_i(t)\right| +\left(1-\rho \underline{a}_{ii}-h_{ii} \right) \mu_i \left|e^1_i(t)\right|  \right\} + \sum \limits_{i=1}^n \sum\limits_{j=1}^m \mu_i\rho \widetilde{a^{\ast}}_{ij} \left|e_j^2(t)\right|  \notag \\
&=& \sum \limits_{i=1}^n \left\{\sum_{j=1,j\neq i}^n \frac{\mu_j}{\mu_i}\left(\left|h_{ji}\right|+\widetilde{a}_{ji}\right) +1-\rho \underline{a}_{ii}-h_{ii} \right\}\mu_i\left|e^1_i(t)\right| + \sum \limits_{j=1}^m \sum\limits_{i=1}^n \frac{\mu_i}{\tau_j}\rho\widetilde{a^\ast}_{ij} \tau_j\left|e_j^{2}(t)\right|. \label{DV2}
\end{eqnarray}
Now from (\ref{DV1}) and (\ref{DV2}), one has
\begin{eqnarray}
&&\leftidx{_0^C}D{_t^\alpha}V_1\left(t^+,e^1\left(t^+\right)\right) \leq \sum \limits_{i=1}^n \left\{\sum_{j=1,j\neq i}^n \frac{\mu_j}{\mu_i}\left(\left|h_{ji}\right|+\widetilde{a}_{ji}\right) +1-\rho \underline{a}_{ii}-h_{ii} \right\}\mu_i\left|e^1_i(t)\right|  \notag \\
&&+ \sum \limits_{j=1}^m \sum\limits_{i=1}^n \frac{\mu_i}{\tau_j}\rho\widetilde{a^\ast}_{ij} \tau_j\left|e_j^{2}(t)\right| + \sum \limits_{i=1}^n \mu_i\left(\left|h_{ii}\right|-1\right)\left|e^1_i(t)\right| \notag \\
&=& \sum \limits_{i=1}^n \left\{\sum_{j=1,j\neq i}^n \frac{\mu_j}{\mu_i}\left(\left|h_{ji}\right| + \widetilde{a}_{ji}\right) +\left|h_{ii}\right|-\rho \underline{a}_{ii}-h_{ii} \right\}\mu_i\left|e^1_i(t)\right|  \notag \\
&&+ \sum \limits_{j=1}^m \sum\limits_{i=1}^n \frac{\mu_i}{\tau_j}\rho\widetilde{a^\ast}_{ij} \tau_j\left|e_j^{2}(t)\right|.       \label{DV12}
\end{eqnarray}
Let
\begin{equation}
V_2\left(t,e^2(t)\right)= \sum \limits_{j=1}^m \tau_j\left|e_j^{2}(t)\right|.  \label{Vt2}
\end{equation}
Then as same as the proof of (\ref{DV12}), by Lemmas \ref{nonexpan} and \ref{Dsgn} with assumptions (A$_3$) and (A$_4$), it follows from (\ref{Dei2}) and (\ref{Vt2}) that
\begin{eqnarray}
&&\leftidx{_0^C}D{_t^\alpha}V_2\left(t^+,e^2\left(t^+\right)\right) \leq  \sum \limits_{j=1}^m \left\{\sum_{i=1,i\neq j}^m \frac{\tau_i}{\tau_j}\left(\left|l_{ij}\right| + \widetilde{b}_{ij}\right) +\left|l_{jj}\right|-\lambda \underline{b}_{jj}-l_{jj} \right\}\tau_j\left|e^2_j(t)\right|  \notag \\
&&+ \sum \limits_{i=1}^n \sum\limits_{j=1}^m \frac{\tau_j}{\mu_i}\lambda\widetilde{b^\ast}_{ji} \mu_i\left|e_i^1(t)\right|.       \label{DV22}
\end{eqnarray}
In light of (\ref{Vt}), (\ref{DV12}) and (\ref{DV22}), one has
\begin{eqnarray}
&&\leftidx{_0^C}D{_t^\alpha}V\left(t^+,\left(e^1\left(t^+\right),e^2\left(t^+\right)\right)\right) = \leftidx{_0^C}D{_t^\alpha}V_1\left(t^+,e^1\left(t^+\right)\right) +\leftidx{_0^C}D{_t^\alpha}V_2\left(t^+,e^2\left(t^+\right)\right)  \notag \\
& \leq &  \sum \limits_{i=1}^n \left\{\sum_{j=1,j\neq i}^n \frac{\mu_j}{\mu_i}\left(\left|h_{ji}\right| + \widetilde{a}_{ji}\right) +\left|h_{ii}\right|-\rho \underline{a}_{ii}-h_{ii} \right\}\mu_i\left|e^1_i(t)\right| \notag \\
&&  \sum \limits_{j=1}^m \left\{\sum_{i=1,i\neq j}^m \frac{\tau_i}{\tau_j}\left(\left|l_{ij}\right| + \widetilde{b}_{ij}\right) +\left|l_{jj}\right|-\lambda \underline{b}_{jj}-l_{jj} \right\}\tau_j\left|e^2_j(t)\right| \notag \\
&& + \sum \limits_{j=1}^m \sum\limits_{i=1}^n \frac{\mu_i}{\tau_j}\rho\widetilde{a^\ast}_{ij} \tau_j\left|e_j^2(t)\right|+ \sum \limits_{i=1}^n \sum\limits_{j=1}^m \frac{\tau_j}{\mu_i}\lambda\widetilde{b^\ast}_{ji} \mu_i\left|e_i^1(t)\right| \notag \\
& = & \sum \limits_{i=1}^n \left\{\sum_{j=1,j\neq i}^n \frac{\mu_j}{\mu_i}\left(\left|h_{ji}\right| + \widetilde{a}_{ji}\right) +\sum\limits_{j=1}^m \frac{\tau_j}{\mu_i}\lambda\widetilde{b^\ast}_{ji} +\left|h_{ii}\right|-\rho \underline{a}_{ii}-h_{ii} \right\}\mu_i\left|e^1_i(t)\right|  \notag \\
&& + \sum \limits_{j=1}^m \left\{\sum_{i=1,i\neq j}^m \frac{\tau_i}{\tau_j}\left(\left|l_{ij}\right| + \widetilde{b}_{ij}\right) +\sum\limits_{i=1}^n \frac{\mu_i}{\tau_j}\rho\widetilde{a^\ast}_{ij} +\left|l_{jj}\right|-\lambda \underline{b}_{jj}-l_{jj} \right\}\tau_j\left|e^2_j(t)\right|  \notag \\
& = & -\left\{ \sum \limits_{i=1}^n \left(1-\xi_i\right)\mu_i\left|e^1_i(t)\right|  + \sum\limits_{j=1}^m \left(1-\zeta_j\right)\tau_j\left|e^2_j(t)\right| \right\}  \notag \\
& \leq & - \theta \left(\sum \limits_{i=1}^{n} \mu_{i}|e_{i}^{1}(t)| + \sum\limits_{j=1}^{m} \tau_{j}|e_{j}^{2}(t)|\right)  \notag \\
& = & - \theta V\left(t,\left(e^1(t),e^2(t)\right)\right),  \notag
\end{eqnarray}
where $\xi_{i}$ and $\zeta_{j}$ are defined by (\ref{xi}) and (\ref{zeta}), respectively, and
\begin{equation*}
\theta =\min \left\{\min\limits_{1\leq i\leq n}\left(1-\xi_{i}\right),\min\limits_{1\leq j\leq m} \left(1-\zeta_{j}\right)\right\}.
\end{equation*}
According to Lemma \ref{MLS}, as same as the proof of Theorem 3 in \cite{ZYW}, we can prove that (\ref{FIIPNN}) is globally Mittag-Leffler stable for each $A \in A_I$, $A^\ast \in A^\ast_I$, $B \in B_I$ and $B^\ast \in B^\ast_I$, i.e.,
$$\left\|\left(x(t),y(t)\right)-\left(x^\ast,y^\ast\right)\right\|\leq V(0,\left(x(0),y(0)\right)-\left(x^\ast,y^\ast\right))E_\alpha\left(-\theta t^\alpha\right),$$
where $\left(x^\ast,y^\ast\right)$ is an equilibrium point of (\ref{FIIPNN}), which completes the proof. \hfill $\square$
\end{theorem}

\begin{corollary}\label{cor3.1}
Assume that
\begin{itemize}
\item[\mbox{(H1)}] $1 - \rho \overline{a}_{ii} \geq 0,\ i=1,2,\cdots,n$;
\item[\mbox{(H2)}] $1 - \lambda \overline{b}_{jj} \geq 0,\ j=1,2, \cdots, m$;
\item[\mbox{(H3)}]  There exist constants $\mu_i>0 \; (i=1,2,\ldots,n)$ and $\tau_j>0 \; (j=1,2,\ldots,m)$ such that
$$\left \{
\begin{array}{l}
0 < 1- \rho \underline{a}_{ii}+ \sum \limits_{j=1,j\neq i}^{n}\frac{\mu _{j}}{\mu _{i}} \rho\widetilde{a}_{ji}+ \sum \limits_{j=1}^{m}\frac{\tau_{j}}{\mu_{i}} \lambda \widetilde{b^{\ast}}_{ji} < 1 , \ i=1,2,\cdots,n, \\
0 < 1- \lambda \underline{b}_{jj}+ \sum \limits_{i=1,i\neq j}^{m}\frac{\tau_{i}}{\tau _{j}} \lambda \widetilde{b}_{ij}+ \sum \limits_{i=1}^{n}\frac{\mu_{i}}{\tau _{j}} \rho \widetilde{a^{\ast}}_{ij} < 1 , \ j=1,2, \cdots, m,
\end{array}
\right.
$$
where
$$
\widetilde{a}_{ji}=\max\left\{\left| \underline{a}_{ji}\right|,  \left|\overline{a}_{ji}\right|\right\}, \quad \widetilde{a^{\ast}}_{ij}=\max\left\{\left|\underline{a^{\ast}}_{ij}\right|, \left|\overline{a^{\ast}}_{ij}\right|\right\}
$$
and
$$
\widetilde{b}_{ij}=\max\left\{\left|\underline{b}_{ij}\right|,\left|\overline{b}_{ij}\right|\right\}, \quad
\widetilde{b^{\ast}}_{ji}=\max\left\{\left|\underline{b^{\ast}}_{ji}\right|, \left|\overline{b^{\ast}}_{ji}\right|\right\}.
$$
\end{itemize}
Then, the system of interval projection neural networks with fractional-order \eqref{e1.2} is globally Mittag-Leffler stable for each $A \in A_I$, $A^\ast \in A^\ast_I$, $B \in B_I$ and $B^\ast \in B^\ast_I$.
\end{corollary}
\begin{remark}
We note that Corollary \ref{cor3.1} is an improved version of Theorem 4.1 in \cite{WZH2}.
\end{remark}

\begin{corollary}\label{cor3.3}
Suppose that
$$1-\sum_{j=1,j \neq i}^n\frac{\mu_j}{\mu_i} \rho\left|a_{ji}\right|-\left|1-\rho a_{ii}\right|>0,\quad \mu_i>0, \, i=1,2,\ldots,n. $$
Then, the global fractional-order projective dynamical system (\ref{e1.3}) is globally Mittag-Leffler stable.
\end{corollary}

\begin{remark}
It is worth mentioning that Corollary \ref{cor3.3} is an improved version of Theorem 4.1 (a) in \cite{WZ}.
\end{remark}

\begin{remark}
Following the paper \cite{YHJ}, we studied the $\alpha$-exponential stability for the global fraction-order projective dynamical system in \cite{WZ} and for the system of fractional-order interval projection neural networks in \cite{WZH2} without noticing the paper \cite{LPG},  in which the authors pointed out that the conclusion of $\alpha$-exponential stability introduced in \cite{YHJ} is invalid in the fractional system. Referring the papers \cite{CZJ,CC,LPG,ZYW}, we would like to point out that the conclusions of $\alpha$-exponential stability in \cite{WZ,WZH2} should be replaced by the Mittag-Leffler stability.
\end{remark}

\section{Numerical examples}
\noindent \setcounter{equation}{0}

This section gives two examples to demonstrate the main results presented in Section 3.

\begin{example}
Assume that $\alpha=0.8$, $\rho=0.3$, $\lambda=0.2$, $a=(-7.1,4.2,-2.4)^\top$, $b=(-3.5, 1.2)^\top$,
\begin{equation*}
\underline{A}=\left(
\begin{array}{ccc}
2.6 & 0.3 & -0.3 \\
-0.5 & 3.4 & -0.1 \\
0.2 & 0.6 & 2.1
\end{array}
\right), \quad \overline{A}=\left(
\begin{array}{ccc}
2.9 & 0.5 & 0.3 \\
-0.4 & 3.6 & 0.2 \\
0.4 & 0.8 & 2.5
\end{array}
\right) ,\quad \underline{A^\ast}=\left(
\begin{array}{cc}
-0.3 & 0.2 \\
0.1 & -0.4 \\
-0.2 & 0.1
\end{array}
\right) ,
\end{equation*}
\begin{equation*}
\overline{A^\ast}=\left(
\begin{array}{cc}
0.2 & 0.4 \\
0.3 & -0.3 \\
0.1 & 0.3
\end{array}
\right),\quad
\underline{B}=\left(
\begin{array}{cc}
3.5 & 0.4 \\
-0.2 & 2.6
\end{array}
\right), \quad \overline{B}=\left(
\begin{array}{cc}
3.6 & 0.7 \\
0.2 & 2.8
\end{array}
\right) ,
\end{equation*}
\begin{equation*}
\underline{B^\ast}=\left(
\begin{array}{ccc}
-0.4 & 0.1 & -0.3 \\
0.5 & -0.2 & 0.6
\end{array}
\right), \quad \overline{B^{\ast }}=\left(
\begin{array}{ccc}
0.5 & 0.3 & 0.4 \\
0.7 & -0.3 & 0.7
\end{array}
\right),
\end{equation*}
\begin{equation*}
\left(h_{ij}\right)_{3\times3}=\left(
\begin{array}{ccc}
0.09 & 0.06 & -0.03 \\
-0.05 & -0.17 & 0.08 \\
0.07 & -0.06 & 0.11
\end{array}
\right), \quad
\left(l_{ij}\right)_{2\times2}=\left(
\begin{array}{cc}
-0.11 & -0.03 \\
-0.08 & 0.09
\end{array}
\right) ,
\end{equation*}
$$K_1(x)=u^1(x)+K_1, \quad u^1(x)=(u^1_1(x),u^1_2(x),u^1_3(x))^\top = \left(h_{ij}\right)_{3\times3}\cdot(x_1,x_2,x_3)^\top,$$
$$K_2(y)=u^2(y)+K_2, \quad u^2(y)=(u^2_1(y),u^2_2(y))^\top = \left(l_{ij}\right)_{2\times2}\cdot(y_1,y_2)^\top,$$
$$K_1=\left\{x \in R^3 | 3 \leq x_1 \leq 4,\, -1.5 \leq x_2 \leq -0.5,\, 0.5 \leq x_3 \leq 1.5 \right\},$$
and
$$K_2=\left\{y \in R^2 | 1.5 \leq y_1 \leq 2.5,\, -2.5 \leq y_2 \leq -1 \right\}.$$

When $\mu_i=\tau_j=1$ $(i=1,2,3, \, j=1,2)$, we know that all the assumptions (A$_1$)-(A$_4$) are satisfied. Therefore, it follows from Theorems \ref{Exist} and \ref{Stab} that  FIIPNN (\ref{FIIPNN}) has a unique equilibrium point for each $A \in A_I$, $A^\ast \in A^\ast_I$, $B \in B_I$ and $B^\ast \in B^\ast_I$ and (\ref{FIIPNN}) is globally Mittag-Leffler stable for each $A \in A_I$, $A^\ast \in A^\ast_I$, $B \in B_I$ and $B^\ast \in B^\ast_I$. Figure \ref{ex1} shows the trajectories of (\ref{FIIPNN}) with the same initial value $x_0=(8.6,-7.3,-5.2)^\top$, $y_{0}=(6.7,-8.5)^\top$ when $A=\underline{A}$, $A^\ast=\underline{A^\ast}$, $B=\underline{B}$, $B^\ast=\underline{B^\ast}$ and $A=\overline{A}$, $A^\ast=\overline{A^\ast}$, $B=\overline{B}$, $B^\ast=\overline{B^\ast}$, respectively.
\begin{figure}[H]
\begin{center}
%Requires \usepackage{graphicx}
\includegraphics[width=3.2in]{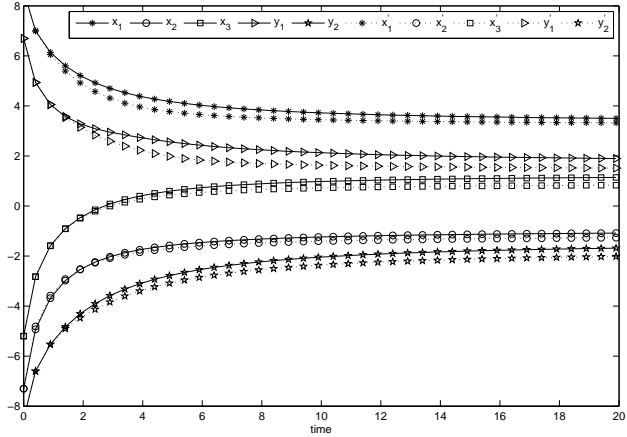}\\[0pt]
\end{center}
\caption{The line $\left(\left(x_1,x_2,x_3 \right)^\top,\left(y_1,y_2\right)^\top\right)$ denotes the transient behavior of FIIPNN (\ref{FIIPNN}) when $A=\protect\underline{A}$, $A^{\ast }=\protect\underline{A^{\ast }}$, $B=\protect\underline{B}$, $B^{\ast }=\protect\underline{B^{\ast }}$. The line $\left(\left(x_1',x_2',x_3'\right)^\top,\left(y_1',y_2'\right)^\top\right)$ denotes the transient behavior of FIIPNN (\ref{FIIPNN}) when $A=\overline{A}$, $A^{\ast}=\overline{A^{\ast}}$, $B=\overline{B}$, $B^{\ast }=\overline{B^{\ast }}$.}
\label{ex1}
\end{figure}
\end{example}

\begin{example}
Let us consider the following fractional-order interval implicit projection neural network
\begin{equation}
\left \{
\begin{array}{l}
\leftidx{_0^C}D{_t^\alpha}x(t)=P_{u(x(t))+K_1}\left[ x(t) -\rho Ax(t)- \rho a \right] -x(t),\quad t\ge 0,\\
x_{i}\left( 0\right) =x_{i0}, \quad i=1,2,
\end{array}
\right.   \label{example2}
\end{equation}
where $x(t)=(x_1(t),x_2(t))^\top$, $\alpha=0.9$, $\rho=0.25$, $u(x(t))=(u_1(x(t)),u_2(x(t)))^\top = \left(h_{ij}\right)_{2\times2} \cdot (x_1(t),x_2(t))^\top$, $a=(-4.8,0)$, $A \in A_I$,
\begin{equation*}
\underline{A}=\left(
\begin{array}{cc}
3.7 & -1.1 \\
-1.8 & 3.1
\end{array}
\right), \quad \overline{A}=\left(
\begin{array}{cc}
4.6 & 1.3 \\
3.8 & 3.4
\end{array}
\right), \quad
\left(h_{ij}\right)_{2\times2}=\left(
\begin{array}{cc}
-0.2 & 0 \\
0 & 0.11
\end{array}
\right),
\end{equation*}
and
$$K_1=\left\{x=(x_1,x_2)^\top \in R^2 | 0 \leq x_1 \leq 2.5,\; 0 \leq x_2 \leq 0.5 \right\}.$$

Clearly, if $\mu_1=2$ and $\mu_2=1$, then
\begin{equation*}
\left \{
\begin{array}{l}
1 - \rho \overline{a}_{11} - h_{11}=0.05>0 ,\\
1 - \rho \overline{a}_{22} - h_{22}=0.04>0 ,\\
0<1- \rho \underline{a}_{11}-h_{11}+|h_{11}|+\frac{\mu_2}{\mu_1}\max\left\{\left|\rho \underline{a}_{21}\right|,  \left|\rho\overline{a}_{21}\right|\right\}=0.95 <1, \\
0<1- \rho \underline{a}_{22}-h_{22}+|h_{22}|+\frac{\mu_1}{\mu_2}\max\left\{\left|\rho \underline{a}_{12}\right|,  \left|\rho\overline{a}_{12}\right|\right\}=0.875<1.
\end{array}
\right.
\end{equation*}
This implies that all the assumptions (A$_1$)-(A$_4$) are satisfied.  Thus, by Theorems \ref{Exist} and \ref{Stab}, we know that neural network (\ref{example2}) has a unique equilibrium point for each $A \in A_I$ and (\ref{example2}) is globally Mittag-Leffler stable for each $A \in A_I$. Figure \ref{ex2} shows the trajectories of (\ref{example2}) with the same initial value $x_0=(5.8,-4.2)^\top$, when $A=\underline{A}$ and $A=\overline{A}$, respectively.
\begin{figure}[H]
\begin{center}
%Requires \usepackage{graphicx}
\includegraphics[width=3.2in]{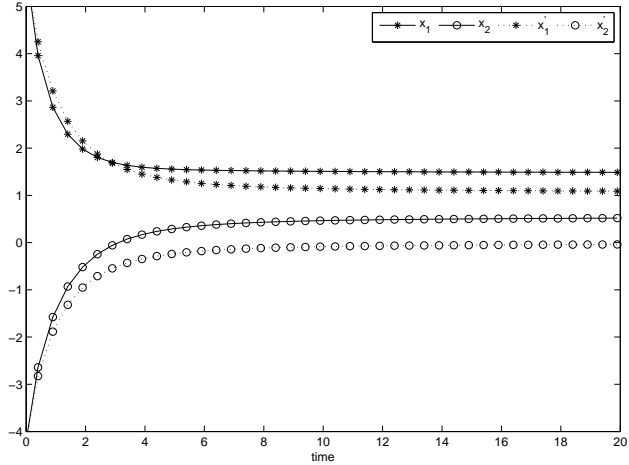}\\[0pt]
\end{center}
\caption{The line $\left(x_1,x_2 \right)^\top$ denotes the transient behavior of (\ref{example2}) when $A=\protect\underline{A}$. The line $\left(x_1',x_2' \right)^\top$ denotes the transient behavior of (\ref{example2}) when $A=\overline{A}$.}
\label{ex2}
\end{figure}
\end{example}


\vskip0.5cm

\begin{thebibliography}{99}
\bibitem{WZH2} Z.B. Wu, Y.Z. Zou, N.J. Huang, A system of fractional-order interval projection neural networks, \textit{J. Comput. Appl. Math.} \textbf{294} (2016) 389-402.

\bibitem{WZ} Z.B. Wu, Y.Z. Zou, Global fraction-order projective dynamical systems, \textit{Commun. Nonlinear Sci. Numer. Simul.} \textbf{19} (2014) 2811-2819.

\bibitem{NNK} M.A. Noor, K.I. Noor, A.G. Khan, Dynamical systems for quasi variational inequalities, \textit{Ann. Funct. Anal.} \textbf{6} (2015) 193-209.

\bibitem{ASA} A.S. Antipin, Minimization of convex functions on convex sets by means of differential equations, \textit{Differ. Equ.} \textbf{30} (1994) 1365-1375.

\bibitem{JB} J. Bolte, Continuous gradient projection method in Hilbert spaces, \textit{J. Optim. Theory Appl.} \textbf{119} (2003) 235-259.

\bibitem{D} K. Ding, N.J. Huang, A new class of interval projection neural networks for solving interval quadratic program, \textit{Chaos Solitons Fractals} \textbf{35} (2008) 718-725.

\bibitem{DN} P. Dupuis, A. Nagurney, Dynamical systems and variational inequalities, \textit{Ann. Oper. Res.} \textbf{44} (1993) 9-42.

\bibitem{FBMT} T.L. Friesz, D. Bernstein, N.J. Mehta, R.L. Tobin, S. Ganjalizadeh, Day-to-day dynamic network disequilibria and idealized traveler information systems, \textit{Oper. Res.} \textbf{42} (1994) 1120-1136.

\bibitem{FSB} T.L. Friesz, Z.G. Suo, D.H. Bernstein, A dynamic disequilibrium interregional commodity flow model, \textit{Transpn Res.-B} \textbf{32} (1998) 467-483.

\bibitem{DG} D. Kinderlehrer, G. Stampacchia, An Introduction to Variational Inequalities and Their Applications, Academic Press, New York, 1980.

\bibitem{WWZ} X.K. Wu, Z.B. Wu, Y.Z. Zou, Existence, uniqueness and stability for a class of interval projective dynamical systems, \textit{Comm. Appl. Nonlinear Anal.} \textbf{20} (2013) 81-94.

\bibitem{WZLX} Z.B. Wu, Y.Z. Zou, X.S. Li, Y.B. Xiao, A new class of global fractional-order projective dynamical systems in Hilbert spaces, \textit{Comm. Appl. Nonlinear Anal.} \textbf{24} (2017) 1-16.

\bibitem{WWZ2} X.K. Wu, Z.B. Wu, Y.Z. Zou, Sensitivity of the set of solutions for a class of fractional set-valued projected dynamical systems, \textit{Nonlinear Anal. Forum} \textbf{20} (2015) 95-105.

\bibitem{WZH1} Z.B. Wu, Y.Z. Zou, N.J. Huang, A class of global fractional-order projective dynamical systems involving set-valued perturbations, \textit{ Appl. Math. Comput.} \textbf{277} (2016) 23-33.

\bibitem{XV} Y.S. Xia, J. Wang, On the stability of globally projected dynamical systems, \textit{J. Optim. Theory Appl.} \textbf{106} (2000) 129-150.

\bibitem{Xia} Y.S. Xia, Further results on global convergence and stability of globally projected dynamical systems, \textit{J. Optim. Theory Appl.} \textbf{122} (2004) 627-649.

\bibitem{Z} D. Zhang, A. Nagurney, On the stability of projected dynamical systems, \textit{J. Optim. Theory Appl.} \textbf{85} (1995) 97-124.

\bibitem{ZS} Y.Z. Zou, C.Y. Sun, Equilibrium points for two related projective dynamical systems, \textit{Comm. Appl. Nonlinear Anal.} \textbf{19}(4) (2012) 109-117.

\bibitem{ZLHS} Y.Z. Zou, X. Li, N.J. Huang, C.Y. Sun, Global dynamical systems involving generalized $f$-projection operators and set-valued perturbation in Banach spaces, \textit{J. Appl. Math.} \textbf{2012} (2012) 12.

\bibitem{ZWZS} Y.Z. Zou, X.K. Wu, W.B. Zhang, C.Y. Sun, An iterative method for a class of generalized global dynamical system involving fuzzy mappings in Hilbert spaces, \textit{Lecture Notes in Commput. Sci.} \textbf{7666}(2012) 44-51.

\bibitem{AHJ} A.A. Kilbas, H.M. Srivastava, J.J. Trujillo, Theory and Applications of Fractional Differential Equations, Elsevier, Amsterdam, 2006.

\bibitem{LZ} C.P. Li, F.R. Zhang, A survey on the stability of fractional differential equations, \textit{Eur. Phys. J. Special Topics} \textbf{193} (2011) 27-47.

\bibitem{OK} N. Ozalp, I. Koca, A fractional order nonlinear dynamical model of interpersonal relationships, \textit{Adv. Difference Equ.} \textbf{2012} (2012) 189.

\bibitem{IP} I. Podlubny, Fractional Differential Equations, Academic Press, San Diego, 1999.

\bibitem{SMM} S.B. Skaar, A.N. Michel, R.K. Miller, Stability of viscoelastic control systems, \textit{IEEE Trans. Automat. Control} \textbf{33} (1988) 348-357.

\bibitem{TB} P.J. Torvik, R.L. Bagley, On the appearance of the fractional derivative in the behavior of real materials, \textit{J. Appl. Mech.} \textbf{51} (1984) 294-298.

\bibitem{Die} K. Diethelm, The Analysis of Fractional Differential Equations: An Application-Oriented Exposition Using Differential Operators of Caputo Type, Springer, Berlin, 2010.

\bibitem{BC} C. Baiocchi, A. Capelo, Variational and Quasi Variational Inequalities, J. Wiley \& Sons, London, 1984.

\bibitem{GLT} R. Glowinski, J.L. Lions, R. Tr\'{e}moli\`{e}res, Numerical Analysis of Variational Inequalities, North-Holland Publishing Company, Amsterdam, 1981.

\bibitem{Mosco} U. Mosco, Implicit Variational Problems and Quasi Variational Inequalities, Nonlinear Operators and the Calculus of Variations, Lecture Notes in Mathematics, Springer-Verlag, New York, 1976.

\bibitem{CP82}D. Chan, J.S. Pang, The generalized quasi-variational inequality problem, \textit{Math. Oper. Res.} \textbf{7} (1982) 211-222.

\bibitem{LCP} Y. Li, Y.Q. Chen, I. Podlubny, Stability of fractional-order nonlinear dynamic systems: Lyapunov direct method and generalized Mittag-Leffler stability, \textit{Comput. Math. Appl.} \textbf{59} (2010) 1810-1821.

\bibitem{ZYW} S. Zhang, Y.G. Yu, H. Wang, Mittag-Leffler stability of fractional-order Hopfield neural networks, \textit{Nonlinear Anal. Hybrid Syst.} \textbf{16} (2015) 104-121.

\bibitem{YHJ} J. Yu, C. Hu, H.J. Jiang, $\alpha$-stability and $\alpha$-synchronization for fractional-order neural networks, \textit{Neural Netw.} \textbf{35} (2012) 82-87.

\bibitem{LPG} K.X. Li, J.G. Peng, J.H. Gao, A comment on ``$\alpha$-stability and $\alpha$-synchronization for fractional-order neural networks'', \textit{Neural Netw.} \textbf{48} (2013) 207-208.

\bibitem{CZJ} J.J. Chen, Z.G. Zeng, P. Jiang, Global Mittag-Leffler stability and synchronization of memristor-based fractional-order neural networks, \textit{Neural Netw.} \textbf{51} (2014) 1-8.

\bibitem{CC} B.S. Chen, J.J. Chen, Global asymptotical $\omega$-periodicity of a fractional-order non-autonomous neural networks, \textit{Neural Netw.} \textbf{68} (2015) 78-88.
\end{thebibliography}
\end{document}